\documentclass[a4paper,11pt]{article}
\usepackage{a4wide}
\usepackage{amssymb}
\usepackage{amssymb}
\usepackage{lipsum}
\usepackage{graphicx}	% Including figure files
\usepackage{amsmath}
\usepackage{subfigure}
\usepackage{multirow}
\usepackage{caption}
\usepackage{soul}
\usepackage{breqn}
\usepackage{tikz,caption}
\usepackage{times}
\usepackage{helvet} \usepackage{courier} 
\begin{document}
\title{Travelling wave solutions of an equation of Harry Dym type arising in the Black-Scholes framework}

\author{Jorge P. Zubelli \\
Department of Mathematics\\ Khalifa University\\ Abu Dhabi, P.O. Box 127788, UAE\\ ADIA Lab\\ Abu Dhabi, P.O. Box 3600, UAE  
\and 
Kuldeep Singh\\
Department of Mathematics\\ Khalifa University\\ Abu Dhabi, P.O. Box 127788, UAE\\
\and 
Vinicius Albani\\ Department of Mathematics\\ Federal University of Santa Catarina\\
Santa Catarina, 88040-900, Brazil
\and Ioannis Kourakis\\ Department of Mathematics\\ Khalifa University\\ Abu Dhabi, P.O. Box 127788, UAE\\
} %%%%%%%%%% % Body of Paper Begins
\date{\today}

\maketitle

\begin{abstract}
The Black-Scholes framework is crucial in pricing a vast number of financial instruments that permeate the complex dynamics of world markets. Associated with this framework, we consider a second-order differential operator $L(x, {\partial_x}) := v^2(x,t) (\partial_x^2 -\partial_x)$ that carries a variable volatility term $v(x,t)$ and which is dependent on the underlying log-price $x$ and a time parameter $t$ motivated by the celebrated Dupire local volatility model. In this context, we ask and answer the question of whether one can find a non-linear evolution equation derived from a zero-curvature condition for a time-dependent deformation of the operator $L$. The result is a variant of the Harry Dym equation for which we can then find a family of travelling wave solutions. This brings in extensive machinery from soliton theory and integrable systems. As a by-product, it opens up the way to the use of coherent structures in financial-market volatility studies.
\\

\\

\\

\noindent
Submitted for publication. 
\end{abstract}

\onecolumn

\section{Introduction}

The iconic Black-Scholes model~\cite{bls} for pricing financial derivative products completed in 2023 its first jubilee with a tremendous impact on world markets~\cite{blanka2023}. In their groundbreaking work, Fisher Black and Myron Scholes proposed a methodology that consisted in producing a replicating portfolio to the financial option that depended on an underlying asset.

Furthermore, they assumed that the stochastic dynamics followed by the underlying asset whose (random) price at time $t$ was given  by $S=S_t$ is described by the stochastic differential equation
\begin{equation}
    dS_t = \mu \, S_t \, dt + \sigma \, S_t \, dW_t \mbox{, }
\end{equation}
where $\mu$ is the so-called drift or expected return on the stock, $\sigma$ is its volatility, and $dW$ is the differential of the Brownian motion $W_t$ in a suitable filtered probability space~\cite{korn}.
Under their methodology, the price of the financial derivative $P=P(S,t)$ satisfied the now celebrated Black-Scholes equation:
\begin{align}
& \frac{\partial P}{\partial t}+
\frac{1}{2}{\sigma}^2S^2\frac{\partial^2 P}{\partial S^2}
+ r \left(S\frac{\partial P}{\partial S}-P\right) = 0 \qquad  \mbox{ for $t\le T_E$ }, \qquad  \mbox{ and }
& P(T_E,\cdot) = h  \mbox{ .}\label{BS0}
\end{align}
Here, $r$ is the risk-free interest rate and $T_E$ is the expiration time of the derivative.
Assume that $\sigma$ and $r$ are constant w.r.t. $t$ and $S$, the final value problem in Eq. (\ref{BS0}) can be reduced to the initial value problem of the heat equation with a suitable change of variables. This assumption allowed for the fast implementation of pricing formulae for the so-called European vanilla options.

However, it turns out that the volatility $\sigma$ is neither constant in time nor at different values of $x$. This phenomenon, known as the implied-volatility {\em smile},  led Bruno Dupire to consider the so-called local volatility models~\cite{dupire}.  In such models, the volatility is subsumed to take the form
$$
\sigma = \sigma (S, t) \mbox{ .}
$$
In this context, the pricing equation becomes
\begin{align} \label{BS1}
& \frac{\partial P}{\partial t}+
\frac{1}{2}{ \sigma(t,S)} ^2 S^2\frac{\partial^2 P}{\partial S^2}
+r \left(S\frac{\partial P}{\partial S}-P\right) = 0  \qquad \mbox{ for } t\le T_E,\\
& P(T_E, S) = h(S) \mbox{ . }
\end{align}
In order to avoid arbitrage opportunities and the mispricing of financial instruments, it is fundamental to calibrate the local volatility consistently with the prices of the so-called calls and puts.
This effort was pursued by a number of authors  starting from Marco~Avellaneda and collaborators~\cite{AFHS97}.

As it turns out $r$ tends to change very little during the life of an option, thus it can be assumed to be constant. Alternatively, one can work with the so-called forward prices. In either case, it is customary in Eq. (\ref{BS1}) to scale out the interest rate $r$ by adopting a change of (the dependent and independent) variables,  so as to concentrate all the efforts on the analysis of the operator
\begin{equation*}
L[P]:= \frac{1}{2} { \sigma(t,S)}^2 \, S^2 \, \frac{\partial^2 P}{\partial S^2} \mbox{ .}
\end{equation*}
%%%%
This operator, after a change of variables to the so-called $log$-prices, becomes
\begin{equation} \label{BS2}
\mathbf{L}[\psi]:= \frac{1}{2}{v^2(x,t) \left(\frac{\partial^2}{\partial x^2} - \frac{\partial}{\partial x}\right ) \psi   } \, ,
\end{equation}
where $v^2(x,t) = \sigma^2(S_0 \exp(x), t)$ and $S = S_0 \exp(x)$.

In this article, we answer the question of whether one can find a nonlinear evolution equation derived from a zero-curvature condition for a time-dependent deformation of the operator $\mathbf{L} $. We shall show that this leads to a variant of the Harry Dym equation \cite{zub2002,zub2008}. We also obtain traveling wave solutions to the corresponding non-linear evolution equation that we have derived.

\section{Preliminaries and Motivation}

As mentioned in the Introduction, a very natural question is how to describe the volatility as a function of time and the value of the underlying in a consistent way with the observed option prices in the market. The importance of this was observed very early and is consistent with the {\em smile effect}; see~\cite{derkani,dermankanizou} and the references therein.
This is the so-called
volatility surface calibration problem~\cite{gatheral2006volatility}.
The literature on the study of the volatility is vast, both from the numerical as well as the theoretical view point. It was initiated very early with the seminal work of Dupire~\cite{dupire} and then with the work of Avellaneda and collaborators~\cite{AFHS97,ABFGKN2000,aveICM, aveSMF, aveIJTAF}; see also~\cite{JacksonSuli98,jermakyan, LagnardoOsher97,stutzer1996simple,samperi}.
Related theoretical aspects were studied by Isakov and collaborators~\cite{BouIsa1997,BouIsa1999,Isa2006} whereas practical aspects were studied by Alex~Lipton {\em et al}~\cite{lipton2002a,lipton2002b,lipton2018,lipton2011}.

As it turns out, the volatility surface calibration problem is an ill-posed inverse problem~\cite{engl1996regularization}. This led to substantial efforts by a number of authors in providing a regularization framework to solve the inverse problem for the volatility surface. These efforts can certainly be traced to the work of Crepey~\cite{Cre2003a,Cre2003b} and subsequently to \cite{achdoupironneau2005,Egger-Engl2005, Hof-Kra-2005,  HKPS-2007} in different contexts. A systematic effort for the calibration of volatility surface from option data using convex regularization techniques has started in the PhD dissertation of A.~De Cezaro~\cite{decezaro2010phd} and has been developed in subsequent work \cite{CezaroScherzeZubelli2012,DZ2015,AlbaniZubelli2014,AlbAscYanZub2015,AlbDecZub2016,AAZ2017,ADZ2017,albani2020splitting} in different contexts.

The nonparametric estimation aspect in the solution of the inverse problem of the volatility considered in the articles \cite{CezaroScherzeZubelli2012,AlbaniZubelli2014,albani2020splitting} naturally leads to the question of finding a good basis for the representation of the volatility. In the present work, we search for building blocks for such a basis among traveling wave solutions of nonlinear evolution equations associated to deformations of the differential operator $\mathbf{L}$ presented in Eq.~\eqref{BS2}. Indeed, coherent structures should play a role in describing the coefficient $\nu$ in the operator  $\mathbf{L}$ at least on a qualitative form.
The main goal of this work is to investigate such coherent travelling waves, which to the best of our knowledge have never been considered in the literature so far, although a very similar operator leads to the Harry Dym equation.

\emph{Solitons}, i.e. exact solutions of integrable PDEs with unique stability properties that allow them to maintain their profile throughtout propagation and mutual interactions (collisions), have played a tremendous role in nonlinear dynamics since the 1960s with the  development of the inverse scattering technique to tackle the Korteweg-de Vries (KdV) equation \cite{drazin1989solitons}. In fact, the first observation of the so-called soliton dates back to 1834 when John Scott Russel, a Scottish engineer, observed a solitary wave travelling down a canal in Scotland without dissipating \cite{dauxois2006physics,remoissenet2013waves}. The fact that equations like the KdV eq. admit stable, localized-profile solutions (solitons) with remarkable stability properties ushered an important chapter in Applied Mathematics.
Further developments to much broader contexts and generalizations such as the Nonlinear Schroedinger (NLS) equation \cite{sulem2007nonlinear}  and the modified KdV equation~\cite{zbMATH03634788} led to a blossoming of a broad area in nonlinear science. Solitonic equations have been found in various fields, from fluid dynamics and optics \cite{dauxois2006physics,remoissenet2013waves}  to quantum field theory and plasma dynamics \cite{mckerr2014freak,singh2022stability,singh2023dust,mushtaq2024nonlinear}. Soliton bearing nonlinear equations have appeared in many unrelated areas, one of which was the so-called bispectral problem posed by Grunbaum and Duistermaat in \cite{duistermaat1986differential} and further developed extensively in other contexts \cite{zubelli1990differential,Zubelli1991,Zubelli1992,AGBrianZubelli}, some of which with deep algebraic connections; see for example~\cite{sym14102202,wilson,VeselovWillox2015} and references therein.

In the soliton literature, the Harry Dym equation appears in a number of studies and has been generalized extensively: see, for instance, in~\cite{zub2002,zub2008} and references therein.

The plan for this article is the following: in the following Section~\ref{mHD}, we derive a modified Harry Dym equation by means of a zero-curvature condition deformation of the operator $\mathbf{L}$ of Eq.~\eqref{BS2}. We shall call this equation the {\em Financial Harry Dym} equation (FHD) due to its motivation from the financial Black-Scholes model. In Section~\ref{Soliton}, we shall obtain  travelling wave solutions to the FHD equation and give a graphical description of its behavior.
We will finally conclude in Section~\ref{conclusion}, proposing some directions for further research.

\section{Derivation of the Financial Harry Dym equation } \label{mHD}
In order to analyze the zero curvature deformations of the operator $\mathbf{L}$ in Equation~\eqref{BS2}, we
consider the ODE
\begin{equation}
\left(\frac{\partial^2}{\partial x^2} - \frac{\partial}{\partial x}\right )\psi=-\frac{\lambda}{v^2} \, \psi\, .\label{der1}
\end{equation}
We follow a similar procedure to the one used extensively in the literature by writing the second order equation as a system in the spirit of \cite{zbMATH03634788}, or the derivation of the HD equation in \cite{hereman1989derivation}.

Here, the spectral parameter $\lambda$ varies in an open domain of the complex plane (or in an interval of the real line) and $v(x,t)$ is a bounded positive function of $x$.
Both $v$ and $\psi$ depend on the parameter $t$. Let $\psi_1=\psi$, $\psi_2=\psi_x$ and $\psi_{2,x}=\psi_{xx}$. Then Eq. (\ref{der1}) can be written as
\begin{equation}
\Psi_x=\mathcal{M}\, \Psi \,  \qquad {\rm where}\qquad \Psi=\begin{bmatrix}
    \psi_1\\
    \psi_2\\
    \end{bmatrix} \quad {\rm and} \qquad \mathcal{M}=\begin{bmatrix}
      0          & 1\\
    -\lambda/v^2 & 1\\
    \end{bmatrix}\label{der2}
\end{equation}
Notice that integrability conditions, which result in solvable nonlinear PDEs, occur when a linear ODE, like the one in Eq. (\ref{der2}), is transformed to
\begin{equation}
\Psi_t=\mathcal{N}\, \Psi \,  \qquad {\rm and} \qquad \mathcal{N}=\begin{bmatrix}
      A   & B\\
      C   & D\\
    \end{bmatrix}\label{der3}
\end{equation}
in a manner that maintains the equation's specific characteristic (such as the spectrum $\lambda$). Consequently, by setting $\lambda_t=0$ (the isospectral case), the integrability conditions: $\Psi_{xt}=\Psi_{tx}$, yield the structure equation
\begin{equation}
 \mathcal{M}_t +[\mathcal{M},\mathcal{N}]=\mathcal{N}_x\, ,   \label{der4}
\end{equation}
where we have employed the commutator $[\mathcal{M},\mathcal{N}]=\mathcal{MN}-\mathcal{NM}$. Using the explicit forms for $\mathcal{M}$ and $\mathcal{N}$ in Eq. (\ref{der4}) leads to
\begin{equation}
\begin{bmatrix}
      0          & 0\\
    -\left(\frac{\lambda}{v^2}\right)_t & 0\\
\end{bmatrix}
+
\begin{bmatrix}
      C+\frac{\lambda}{v^2}B    & D-(A+B)\\
    -\frac{\lambda}{v^2}(A-D)+C   & -C-\frac{\lambda}{v^2}B   \\
\end{bmatrix}
=
\begin{bmatrix}
      A_x    &  B_x\\
      C_x    &  D_x\\
\end{bmatrix}\label{der5}
\end{equation}
This will give four coupled equations, one of which describes the evolution of $r(x, t )$:
\begin{eqnarray}
    C+\frac{\lambda}{v^2}B=A_x  \label{der5a}\\
    D-(A+B)=B_x   \label{der5b}\\
    \frac{2\lambda
     v_t}{v^3}-\frac{\lambda}{v^2}(A-D)+C=C_x  \label{der5c}\\
    -C-\frac{\lambda}{v^2}B=D_x   \, .\label{der5d}
\end{eqnarray}
From Eq. (\ref{der5a}) and (\ref{der5d}), it follows that
\begin{equation}
D=-A \label{der6a} \, .
\end{equation}
Now, Eqs. (\ref{der5a}) and (\ref{der5b}) can be re-written as
\begin{eqnarray}
    C=-\frac{1}{2}(B_{xx}+B_x)-\frac{\lambda}{v^2}B ,  \label{der6b}\\
    A=-\frac{1}{2}(B_x+B)\,
      \, .\label{der6c}
\end{eqnarray}
Using Eqs. (\ref{der6a})-(\ref{der6c}) in Eq. (\ref{der5c}) so that the evolution of $v$ can be written as
\begin{equation}
\frac{v_t}{v^3}=-\frac{B_{xxx}}{4\lambda}+\frac{B_{x}}{4\lambda}+\frac{v_x}{v^3}B-\frac{B_x}{v^2}\, .\label{der7}
\end{equation}

Let $B=-4\lambda v$  leads to a modified Harry-Dym (HD) equation in the form:
\begin{equation}
 \frac{v_t}{v^3}=v_{xxx}-v_x \, .\label{der8}
\end{equation}

\section{The Financial Harry Dym equation and its Travelling-Wave Solution} \label{Soliton}
%%%%

The modified Harry Dym equation derived above reads
\begin{equation}
    v_{t}=v^3 (v_{xxx}-v_{x}) \, . \label{e1}
\end{equation}

Anticipating stationary-profile localized solutions, we may consider a moving reference frame, by introducing the ``traveling coordinate" $\xi=x-\Lambda t$, where $\Lambda$ is a real parameter, whose value is left arbitrary at this stage.
Therefore, the partial derivatives in the above PDE will transform as  $\frac{\partial \cdot}{\partial x}\rightarrow \frac{d \cdot}{d \xi}$ and $\frac{\partial \cdot}{\partial t}\rightarrow -\Lambda\frac{d \cdot}{d \xi}$.

Now, Eq. ({\ref{der8}}) can be written as
\begin{eqnarray}
 -\frac{\Lambda}{v^3}\frac{d v}{d\xi}+\frac{d v}{d\xi}=\frac{d^3 v}{d\xi^3},\\
 \frac{\Lambda}{2}\frac{d }{d\xi}\left(\frac{1}{v^2}\right)+\frac{d v}{d\xi}=\frac{d^3 v}{d\xi^3},\\
\int \frac{d}{d \xi}\left[\frac{\Lambda}{2v^2}+v \right]d\xi=\int \frac{d^3 v}{d\xi^3} d\xi,\\
\frac{\Lambda}{2v^2}+v=\frac{d^2 v}{d\xi^2}+c_1 \, ,\label{e3}
\end{eqnarray}
where $c_1$ is an arbitrary real constant. Using suitable boundary conditions (for localized modes), viz. $v\rightarrow v_0$, $v^{'}\rightarrow 0$ and $v^{''}\rightarrow 0$  as $\xi\rightarrow\pm\infty$, where $v_0$ is a suitable asymptotic value (a constant),  one finds  $c_1=\frac{\Lambda}{2 v_{0}^2}+v_{0}$.

Substituting into Eq.~(\ref{e3}), one can write
\begin{equation}
   \frac{\Lambda}{2}\left(\frac{1}{v^2}-\frac{1}{v_{0}^2}\right)+(v-v_0)=\frac{d^2 v}{d\xi^2} \, ,\label{e4}
\end{equation}

Multiplying by $d v/d\xi$ and then integrating,
\begin{eqnarray}
\int \left(\frac{d^2 v}{d\xi^2}\right)\left(\frac{d v}{d\xi}\right) d\xi=\int  \frac{\Lambda}{2}\left(\frac{1}{v^2}-\frac{1}{v_{0}^2}\right) dv +\int (v-v_0) dv\, ,\\
\frac{1}{2} \left(\frac{d v}{d \xi}\right)^2 +c_2 = -\frac{\Lambda}{2 v}-\frac{\Lambda v}{2v_{0}^2}+\frac{v^2}{2}-v\, v_0 \, , \label{e5}
\end{eqnarray}
where an arbitrary real constant $c_2$ has been defined, upon integration in the latter step.
Again, using vanishing boundary conditions (for localized modes), i.e. assuming that equilibrium holds far from the excitation, viz.
\[ v\rightarrow v_0 \qquad {\rm and} \qquad v^{'}\rightarrow 0 \qquad {\rm as} \qquad  \xi\rightarrow\pm\infty \, ,\] we get  $c_2=-\frac{\Lambda}{v_{0}}-\frac{v_{0}^2}{2}$.
Eq. (\ref{e5}) can thus be written as
\begin{equation}
\frac{1}{2} \left(\frac{d v}{d \xi}\right)^2 = -\frac{1}{2v \, v_0^2}(\Lambda - v\, v_0^2)(v-v_{0})^2 \, .
\end{equation}
This relation can be cast in the form:
\begin{equation}
\frac{1}{2} \left(\frac{d v}{d \xi}\right)^2 +S(v)=0\, .\label{e6}
\end{equation}
where the pseudopotential function $S$ is given by
\begin{equation}
 S(v)=\frac{1}{2 v\, v_0^2}(\Lambda-v\, v_0^2)(v-v_0)^2 . \label{e7}
\end{equation}

\begin{figure*}
\centering
\subfigure[]{\includegraphics[width=2.5in]{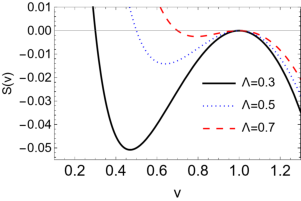}}\hspace{0.8cm}
\subfigure[]{\includegraphics[width=2.5in]{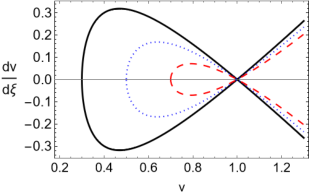}}\\
\subfigure[]{\includegraphics[width=2.5in]{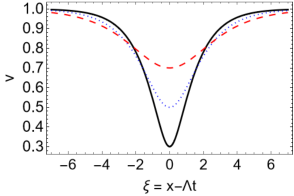}}
\caption{Plots of (a) the pseudopotential function, (b) phase portrait ($v_{\xi}$) and (c) the associated soliton form ($v$), for different values of $\Lambda$. Here, we have taken $v_0 = 1$ in the parametric analysis.}\label{f1}
\end{figure*}

\begin{figure*}
\centering
\subfigure[]{\includegraphics[width=3in]{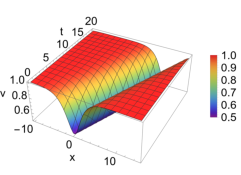}}\hspace{0.8cm}
\subfigure[]{\includegraphics[width=2.5in]{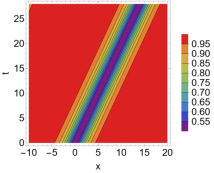}}
\caption{(a) A  3D  plot of v(x,t) is presented, in the $x-t$ plane. (b) A contour plot corresponding to the (numerically obtained) localized solution is shown, for a fixed value of $\Lambda=0.5 $ and $v_0=1$.}\label{f2}
\end{figure*}

Eq. (\ref{e7}) satisfies the conditions for the existence of localized (soliton-like) solutions
i.e.,
\begin{equation}
 (i) \, S(v)|_{v=v_0}=S'(v)|_{v=v_0}=0 \,  \qquad {\rm and} \qquad (ii) \, S^{''}(v)|_{v=v_0}< 0 \, .
 \end{equation}
The first (two) relations are readily satisfied. From the last one, we get $\Lambda < \, v_{0}^3$, which delimits the existence domain for solitons. In Figs. \ref{f1} and \ref{f2}, we have taken $v_0=1$. Therefore, the existence domain for solitons is $\Lambda \in (0, 1)$.

Based on this, we have plotted the Sagdeev pseudopotential profile $S(v)$  versus $v$ for different indicative values of $\Lambda$ in Fig. \ref{f1}(a).
Fig. \ref{f1}(b) depicts the phase portrait, i.e. $v \, '$ versus $v$, for different values of $\Lambda$.
Fig.\ref{f1}(c) depicts the soliton profile, i.e. $v(\xi)$ versus $\xi$  (obtained upon numerical integration of the above ODE),  for different values of $\Lambda$.
It is evident in all figures that the soliton width decreases while its maximum amplitude increased, as $\Lambda$ increases.

Fig. \ref{f2}(a) illustrates the travelling wave solution for the volatility $v(x,t)$ on the $x-t$ plane, for fixed $\Lambda=0.5$. The corresponding contour plot is also given in panel (b), in the same Figure.

\section{Conclusion and Perspectives} \label{conclusion}

\begin{figure*}
\centering
\subfigure[]{\includegraphics[width=3in]{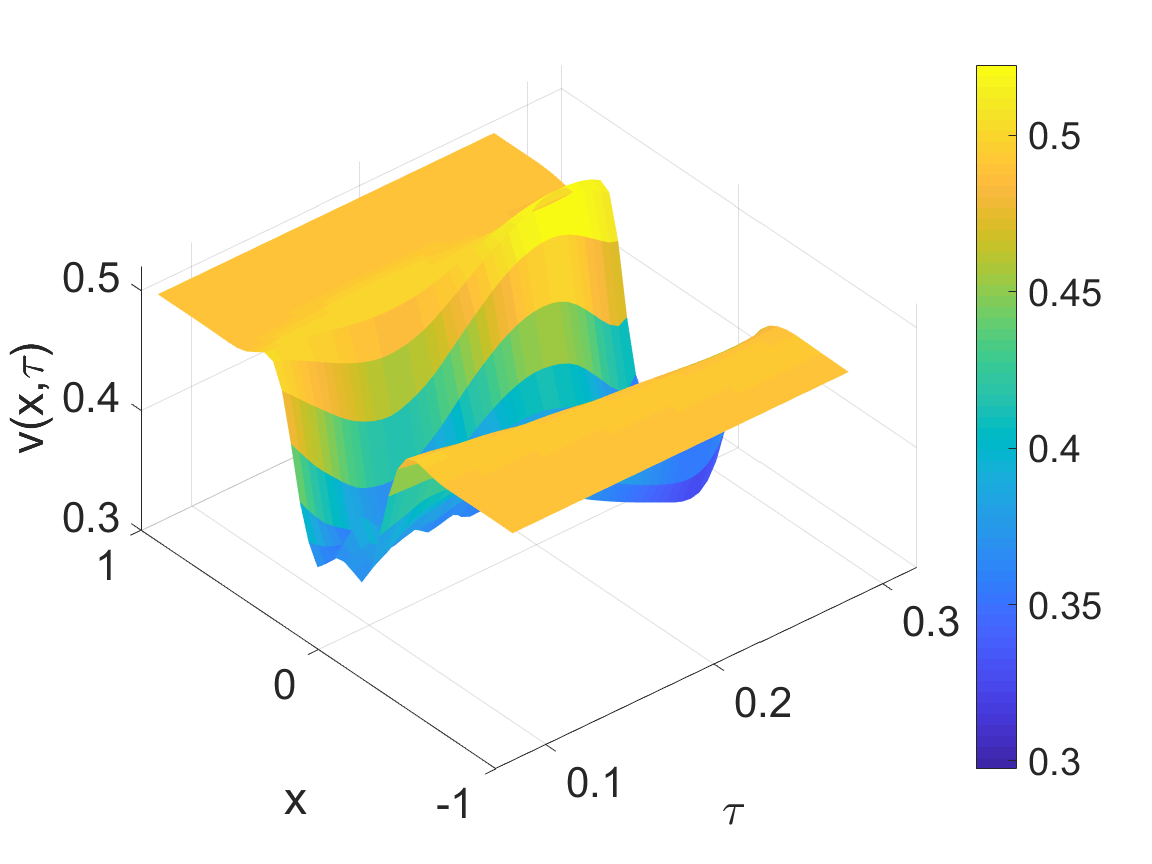}}\hspace{0.8cm}
\subfigure[]{\includegraphics[width=3in]{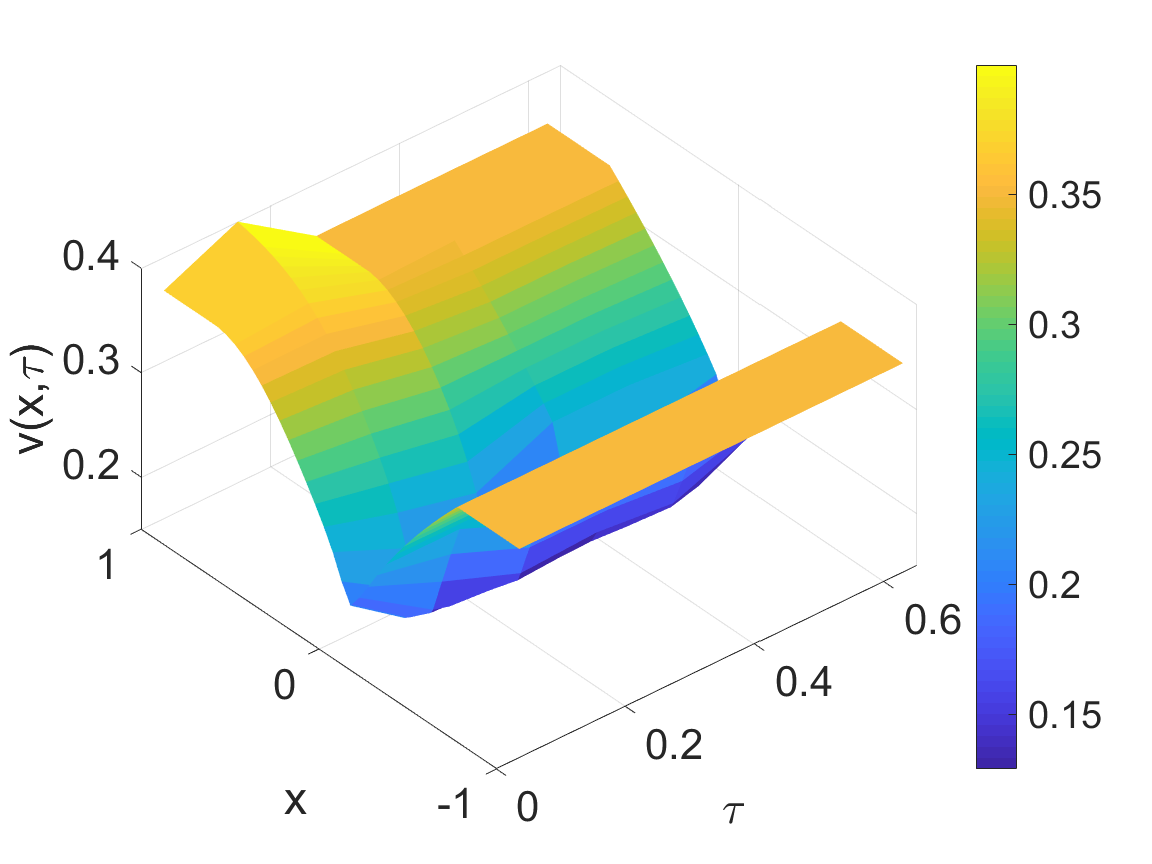}}
\caption{3D  plots of the local volatility surface $v(x,\tau)$ are presented, in the $x-\tau$ plane, where $x$ is the log-moneyness and $\tau$ is the time until expiration. The surfaces in (a) and (b) were reconstructed from call option prices of future contracts of Henry Hub gas and WTI oil, respectively.}\label{f3}
\end{figure*}

The travelling wave solution we obtained for the FHD equation has a remarkable qualitative similarity to the reconstructed volatility surfaces obtained in previous works \cite{AlbDecZub2016,AlbZub2014,AlbAscZub2016,AlbDecZub2017}.
We display two examples of volatility surfaces reconstructed from call options on future contract prices of Henry Hub natural gas and WTI oil in Fig.~\ref{f3}. The surfaces are plotted in the $x-\tau$ plane, where $x$ is the log-moneyness, i.e., $x = \log\left(F/K \right)$ with $K$ the option's strike price and $F$ the underlying future contract price, and $\tau = T_E - t$ is the time to expiration. Henry Hub options were traded on 16-Nov-2011 and WTI options were traded on 04-Sep-2013. The local volatility surfaces were reconstructed using the so-called online Tikhonov-type regularization \cite{AlbZub2014,AlbAscZub2016}, with the regularization parameter chosen according to a discrepancy principle. We remark that since our solution to the inverse problem is not dimensionless, we cannot directly match the FHD solution and the numerical results.

Dupire's equation \cite{dupire} was formulated using the log-moneyness and time-to-expiration variables. The PDE was solved using a Crank-Nicolson scheme where the option prices were normalized by the futures. The log-moneyness mesh ranged from $-5$ to $5$. The time-to-maturity mesh ranged from zero to the maximum maturity time in the sample. The mesh of the observed prices is coarser than the mesh used to solve the PDE. Thus, observed prices were interpolated in the log-moneyness mesh, and the time-to-maturity mesh accounted for the data maturity times. The penalization non-linear least-squares (Tikhonov) functional was minimized using a gradient-descent method. For each commodity, multiple surfaces of prices were used in the estimation, increasing the observed data. The regularity of the families of surfaces of prices and local volatilities with respect to some index was also assumed and included in the estimation through an additive penalty term in the Tikhonov functional \cite{AlbZub2014,AlbAscZub2016}.

In conclusion, the results obtained here suggest that the connection between multi-soliton solutions of nonlinear evolution equations such as the FHD equation should be investigated.

In future work, we plan to explore further the complete integrability of the FHD equation as well as the equations obtained from the deformations of the operator $\mathbf{L}$ in Eq.~\eqref{BS2}. We conjecture that multisoliton solutions could be used as natural basis for the reconstructions presented in Figure~\ref{f3}. However, this would require a careful comparison with the units present in the Black-Scholes model.
\section*{Data availability}

The data sets generated during and/or analyzed during the current study are available from the corresponding author on reasonable request.

\section*{Acknowledgements}
Authors KS and IK gratefully acknowledge financial support from Khalifa University via CIRA (Competitive Internal Research Award) grant CIRA-2021-064/8474000412.

We dedicate this paper to Prof. Pierre Sabatier. It was in one of his papers \cite{Sabatier79} that JPZ first came across the Harry Dym equation.  
JPZ acknowledges a number of encouraging conversations with Dr. Roman Raykov Paunov that stimulated many ideas presented here. 

\medskip

\section*{Author contributions statement}
All authors have contributed equally to the conceptualization, formal design, and methodology of the study. All authors contributed to the analysis of the results. All authors have read and approved the final manuscript.

\section*{Declaration of competing interest}

The authors do not have any kind of conflict of interest.

\bibliographystyle{unsrt}

\bibliography{new2019,solitonsreferences,slvslides,old}

% \begin{thebibliography}{}

% \bibitem[Kawamoto (1985)]{kaw85} Kawamoto S., An exact transformation from the Harry Dym equation to modified KdV equation, \emph{J. Phys. Soc. Japan}, {\bf 54}, 2055 (1985).
% \end{thebibliography}

\end{document}